\documentclass[oneside]{amsart}
\usepackage{enumerate,amssymb}
\usepackage{mathrsfs}

\newtheorem{thm}{Theorem}[section]
\newtheorem{coro}[thm]{Corollary}
\newtheorem{prop}[thm]{Proposition}

\newtheorem{claim}{Claim}

\newtheorem*{thms}{Theorem}
\newtheorem{lem}[thm]{Lemma}

\newtheorem{question}[thm]{Question}
\theoremstyle{definition}
\newtheorem{defn}[thm]{Definition}

\newcommand{\Rset}{\mathbb{R}}
\newcommand{\Nset}{\omega}
\newcommand{\nset}{\omega_+}
\newcommand{\abs}[1]{\lvert#1\rvert}
\newcommand{\Abs}[1]{\lVert#1\rVert}
\newcommand{\seq}[1]{\langle#1\rangle}
\newcommand{\hm}{\mathscr H}
\newcommand{\eps}{\varepsilon}
\newcommand{\del}{\delta}
\newcommand{\subs}{\subseteq}
\renewcommand{\leq}{\leqslant}
\renewcommand{\geq}{\geqslant}
\DeclareMathOperator{\hdim}{\dim_{\mathsf{H}}}
\DeclareMathOperator{\lpdim}{\underline{dim}_{\mathsf{P}}}
\DeclareMathOperator{\lbdim}{\underline{dim}_{\mathsf{B}}}
\DeclareMathOperator{\diam}{diam}
\DeclareMathOperator{\lip}{lip\mkern-1.5mu}
\DeclareMathOperator{\Lip}{Lip\mkern-1.5mu}

\DeclareMathOperator{\por}{\mathsf{por}}

\newenvironment{enum}{\begin{enumerate}[\rm(i)]}{\end{enumerate}}
  {\begin{list}{\textbullet}{\labelwidth1ex\setlength{\leftmargin}{1.5em}}}%
  {\end{list}}
\newenvironment{cproof}%
  {\noindent$\boldsymbol\vdash$}{\renewcommand{\qedsymbol}%
  {$\boldsymbol\dashv$}\hfill\qed\par}
\newcommand{\si}{$\sigma$\nobreakdash-}
\newcommand{\upto}{{\nearrow}}
\newcommand{\downto}{\downarrow}
\newcommand{\VV}{\mathcal V}
\newcommand{\EE}{\mathcal{E}}
\newcommand{\GG}{\mathscr{G}}

\newcommand{\II}[1]{[0,1]^#1}

\newcommand{\Class}{\mathscr C}
\newcommand{\Flass}{\mathscr F}
\newcommand{\class}{\boldsymbol{\mathcal C}}
\newcommand{\klass}{\boldsymbol{\mathcal{K}}}
\newcommand{\gauge}{\psi}
\newcommand{\nuu}{\boldsymbol{\nu}}
\newcommand{\boxm}{\nuu}
\newcommand{\lboxm}{\boxm}

\newcommand{\mc}{\mathcal}
\newcommand{\cross}[1]{\boldsymbol\times\mkern-2mu#1}
\newcommand{\stp}{strongly porous}
\newcommand{\sistp}{\si strongly porous}

\begin{document}
\title
[Little lip of a typical function]
{Little lip of a typical function\\ and \si porosity of graphs}
\author{Ond\v rej Zindulka}
\address
{Ond\v rej Zindulka\\
Department of Mathematics\\
Faculty of Civil Engineering\\
Czech Technical University\\
Th\'akurova 7\\
160 00 Prague 6\\
Czech Republic}
\email{ondrej.zindulka@cvut.cz}
\urladdr{http://mat.fsv.cvut.cz/zindulka}

\begin{abstract}
For a mapping $f\colon X\to Y$ between metric spaces
the function $\lip f\colon X\to[0,\infty]$ defined by
\[
  \lip f(x)=\liminf_{r\to0}\frac{\diam f(B(x,r))}{r}
\]
is termed the \emph{little lip} function of $f$.
We prove that, given a locally compact set
$\Omega\subs\Rset^d$ with a nonempty interior, for a typical continuous function
$f\colon \Omega\to\Rset$ the set $\{x\in\Omega:\lip f(x)>0\}$
has both Hausdorff and lower packing dimensions exactly $d-1$,
while the set $\{x\in\Omega:\lip f(x)=\infty\}$
has non-\si finite $(d{-}1)$-dimensional Hausdorff measure.
This sharp result roofs previous results of
Balogh and Cs\"ornyei~\cite{MR2213746}, Hanson~\cite{MR2929022} and
Buczolich, Hanson, Rmoutil and Z\"urcher~\cite{MR3984283}.
It follows, e.g., that a graph of a typical function $f\in C(\Omega)$ is
\si strongly porous,
microscopic, and for a typical function $f\colon[0,1]\to[0,1]$
there are sets $A,B\subs[0,1]$ of lower packing and Hausdorff dimension
zero, respectively, such that the graph of $f$ is contained in the set
$A\times[0,1]\cup[0,1]\times B$.
\end{abstract}

\maketitle

\section{Introduction}

For a mapping $f\colon X\to Y$ between metric spaces
the function $\lip f\colon X\to[0,\infty]$ defined by
\begin{equation}\label{ll1}
  \lip f(x)=\liminf_{r\to0}\frac{\diam f(B(x,r))}{r}
\end{equation}
is termed the or \emph{little lip} function of $f$.
($B(x,r)$ denotes the closed ball of radius $r$ centered at $x$.)
Let us mention that some authors, cf.~\cite{MR2213746}, define the scaled
oscillations from the version of $\omega_f$ (cf.~\eqref{osc}) given by
$\widehat\omega_f(x,r)=\sup_{y\in B(x,r)} \abs{f(y)-f(x)}$
that may be more suitable for, e.g., differentiation considerations.
It is clear though that
$\widehat\omega_f(x,r)\leq\omega_f(x,r)\leq2\widehat\omega_f(x,r)$,
and thus the two lower scaled oscillation functions differ at most by a factor of $2$.
Since we are interested only in the sets
$\{\lip f=0\}$, $\{\lip f>0\}$ and $\{\lip f=\infty\}$,
it does not matter which of the two definitions we use.

We are interested in the behavior of $\lip f$ for continuous functions
on a locally compact set in the $d$-dimensional Euclidean space.

There is also the \emph{big Lip function} of $f$
\[
  \Lip f(x)=\limsup_{r\to0}\frac{\diam f(B(x,r))}{r}
\]
that plays a crucial role in the
famous Radema\-cher--Stepanov Theorem that asserts that if $f$ is measurable, then
it is differentiable at almost every point where $\Lip f$ is finite.
The little lip function was used in~\cite{MR1708448}
and ~\cite{MR2041901} in the same context: in the study of validity of
the Rademacher Theorem in metric spaces.

Recently there has been a lot of interest in the behavior of functions
at points where $\lip f(x)$ is finite and in particular in functions
that have only a few points with $\lip f(x)=\infty$.
In particular, the differentiability of such functions has been studied and
the structure of exceptional sets $\{x:\lip f(x)=\infty\}$ has been investigated.
Balogh and Cs\"ornyei brought deep results on $\lip f$ in~\cite{MR2213746}
and ignited thus more interest in the notion.
The subject was further studied in~\cite{MR3511937,MR4618645,MR3918835,MR3984283,MR2929022,MR3984283,Hanson2018,Hanson2025,MR4093796}.

The papers~\cite{MR2213746,MR3984283,MR2929022} are of particular importance for
the present paper. In~\cite{MR2213746}, Balogh and Cs\"ornyei constructed
an example contrasting the Rademacher--Stepanov Theorem:
\begin{thms}[{\cite[Theorem 1.4]{MR2213746}}]
There is a continuous function $f\colon[0,1]\to\Rset$ that has $\lip f=0$ almost everywhere
and yet is nowhere differentiable.
\end{thms}
This result was improved in two directions. Hanson's~\cite{MR2929022}
example has even a smaller exceptional set:
\begin{thms}[{\cite[Theorem 2.3]{MR2929022}}]
There is a continuous function $f\colon [0,1]\to\Rset$ that has $\lip f=0$ everywhere
except on a set of Hausdorff dimension zero
and yet it is nowhere differentiable.
\end{thms}
In another direction, Buczolich, Hanson, Rmoutil and Z\"urcher proved in~\cite{MR3984283}
that the example of Balogh and Cs\"ornyei is actually typical
(which sharply contrasts the classical result of Banach~\cite{Ban31}
that a typical continuous function has $\Lip f=\infty$ at each point).
Recall that a property of continuous functions is called \emph{typical} (or \emph{generic}) if
it holds for all functions in a dense $G_\delta$-set in $C([0,1])$.
\begin{thms}[{\cite[Theorem 4.2]{MR3984283}}]
A typical continuous function $f\colon [0,1]\to\Rset$ has $\lip f=0$ almost everywhere.
\end{thms}
Thus a typical function has $\lip f=0$ a.e.~\emph{and} is nowhere differentiable.
Both of the latter results have also higher dimensional versions.

\subsection*{Typical exceptional sets}

The primary goal of this paper is
to prove a sharp theorem roofing the above three results.
In its simplest form it reads as follows. (The full strength
version is presented in Section~\ref{sec:size}).
\begin{thm}\label{thm1.1}
For a typical continuous function $f\colon[0,1]^d\to\Rset$
the exceptional sets
$$
  \{x\in\Omega:\lip f(x)>0\}, \quad \{x\in\Omega:\lip f(x)=\infty\}
$$
are \sistp{} and their lower packing and Hausdorff dimensions are
$d-1$.
\end{thm}
This theorem is proved in Section~\ref{sec:size}.

\subsection*{Level sets}
Maybe surprisingly, information about level sets of a typical function can be derived
from this behavior of $\lip f$.
The study of the size of level sets of a typical function was initiated by
Krishna M.~Garg~\cite{MR143855} and since then it has been extensively studied by many.
To name a few,
Bernd Kirchheim~\cite{MR1260171} proved that all level sets of a typical function
$f\in C([0,1]^d)$ have Hausdorff dimension $d-1$. This result
was recently extended to general compact metric spaces by
Balka, Buczolich and Elekes in papers~\cite{MR3318168,MR3000710}.
On the other hand, Brian S.~Thomson~\cite{MR807990} proved that all level sets
of a typical function $f\in C([0,1])$ are \si strongly porous.

Roughly speaking, these theorems say that level sets of a typical function are small.
The following theorem claims that most level sets of a typical function
are actually contained in a \emph{single} small set. See Section~\ref{sec:level}.
\begin{thm}
For a typical continuous function $f\colon\Omega\to\Rset$,
where $\Omega\subs\Rset^d$ is a locally compact set,
there is set $A\subs\Omega$ such that
\begin{enum}
\item $A$ is \sistp{} and the lower packing dimension of $A$ is at most $d-1$,
\item $f^{-1}(y)\subs A$ for all $y$ except of a set of Hausdorff dimension zero.
\end{enum}
\end{thm}
\subsection*{Size of a graph}

The above result on level sets can be easily rephrased as follows.
\begin{thm}
For a typical continuous function $f\colon\Omega\to\Rset$,
where $\Omega\subs\Rset^d$ is a locally compact set,
there are sets $A\subs\Omega$ and $B\subs\Rset$ such that
\begin{enum}
\item $A$ is \sistp{} and the lower packing dimension of $A$ is at most $d-1$,
\item Hausdorff dimension of $B$ is zero,
\item $f\subs (A\times\Rset)\cup(\Omega\times B)$.
\end{enum}
\end{thm}

This sounds particularly strange in dimension 1:
a graph of a typical function on $[0,1]$ is contained in a rather sparse grid
of horizontal and vertical lines.
\begin{coro}
For a typical function $f\in C([0,1])$ there are sets $A,B\subs[0,1]$
such that $\lpdim A=0$, $\hdim B=0$
and such that $f\subs (A\times\Rset)\cup(\Rset\times B)$.
\end{coro}

Maybe the most surprising consequence of Theorem~\ref{thm1.1} is that
a typical graph is \sistp.
\begin{thm}
The graph of a typical continuous function $f\colon\Omega\to\Rset$,
where $\Omega\subs\Rset^d$ is a locally compact set,
is \sistp{} in $\Rset^{d+1}$.
\end{thm}

These facts regarding the size of a typical graph are discussed in
Section~\ref{sec:graph}.
In Section~\ref{sec:micro} we prove that the exceptional sets, level sets and
the graph of a typical function are microscopic, a notion that has
been recently studied in a number of papers,
see Definition~\ref{def:micro}.

\section{Background: fractal measures, dimension, porosity}\label{sec:prelim}

In this section we review some background and preliminary material.
Namely we recall the notions Hausdorff and lower packing dimensions
and corresponding measures.

Some of the common notation includes
$B(x,r)$ for the closed ball centered at $x$, with radius $r$;
$d$ is a generic symbol for a metric;
we write $\diam E$ for the diameter of a set $E$ in a metric space.
Letters $d,n,m,i,j,k$ are generic symbols for positive integers.
$\Rset$ denotes the real line and $\Rset^d$ the $d$-dimensional Euclidean space;
$\Nset$ stands for the set of natural numbers including zero and
$\nset$ stands for the set of natural numbers excluding zero.
$\abs A$ denotes the cardinality (finite or infinite) of the set $A$.
We write $X_n\upto X$ for an increasing sequence of
sets $\seq{X_n}$ with union $X$.

\subsection*{Hausdorff measure}

Besides standard Hausdorff measures we will also make use of those
induced by gauges.
A non-decreasing, right-continuous function $\gauge\colon [0,\infty)\to[0,\infty)$
such that $\gauge(0)=0$ and $\gauge(r)>0$ if $r>0$ is called a \emph{gauge}.
The following is the common (partial) ordering of gauges, cf.~\cite{MR0281862}:
\[
  \phi\preccurlyeq\gauge\ \overset{\mathrm{def}}{\equiv}
  \ \limsup_{r\to0+}\frac{\gauge(r)}{\phi(r)}<\infty,\qquad
  \phi\prec\gauge\ \overset{\mathrm{def}}{\equiv}
  \ \limsup_{r\to0+}\frac{\gauge(r)}{\phi(r)}=0.
\]
In the case when $\gauge(r)=r^s$ for some $s>0$,
we write $\phi\preccurlyeq s$ instead of $\phi\preccurlyeq\gauge$
and likewise for $\phi\succcurlyeq s$.
%

If $\del>0$, a cover $\mc A$ of a set $E\subs X$ is termed a
\emph{$\del$-cover} if $\diam A\leq\del$ for all $A\in\mc A$.
If $\gauge$ is a gauge,
the \emph{$\gauge$-dimensional Hausdorff measure} $\hm^\gauge(E)$ of
a set $E\subs X$ is defined thus:
For each $\del>0$ set
\[
  \hm^\gauge_\delta(E)=
  \inf\left\{\sum\nolimits_n\gauge(\diam E_n):
  \text{$\{E_n\}$ is a countable $\delta$-cover of $E$}\right\}
\]
and put
\[
  \hm^\gauge(E)=\sup_{\delta>0}\hm^\gauge_\delta(E).
\]
In the common case when $\gauge(r)=r^s$ for some $s>0$, we write $\hm^s$ for
$\hm^\gauge$, and the same licence is used for other measures and set functions
arising from gauges.
We also define the singular case of $\hm^0$ to be the counting measure.
Properties of Hausdorff measures are well-known, see, e.g., \cite{MR0281862}.
%
\subsection*{Hausdorff dimension}
%
Recall that the \emph{Hausdorff dimension} of $E$ is denoted and defined by
\[
  \hdim E=\sup\{s:\hm^s(E)>0\}.
\]
Properties of Hausdorff dimension are well-known, see, e.g., \cite{MR2118797}.

\subsection*{Lower box-counting measure}
Besides Hausdorff dimension we will also
make use of the lower packing dimension.
There are several ways to define it. We do so via the so called
lower box-counting measure.

For $E\subs X$ and $\del>0$, define the \emph{box-counting function} of~$E$
\begin{equation}\label{eq:capacity1}
  N_\del(E)=\inf\{\abs\EE:\EE\text{ is $\del$-cover of $E$}\}.
\end{equation}
For purely technical reasons, we introduce an auxiliary notion:
a right-continuous function $\zeta\colon (0,\infty)\to(0,\infty)$
such that $\varliminf_{r\to0}\zeta(r)=0$
is called a \emph{pseudogauge}.

\begin{defn}[{\cite{MR827889,MR3114775}}]
Let $\zeta$ be a pseudogauge and $E\subs X$. For $\eps>0$ let
\begin{align}
  \boxm^\zeta_\eps(E)&=\inf_{\del\leq\eps}N_\del(E)\zeta(\del),\label{hs1}\\
  \boxm^\zeta_0(E)&=\sup_{\eps>0}\boxm^\zeta_\eps(E)
  =\liminf_{\del\to0}N_\del(E)\zeta(\del),\label{hs2}
\end{align}
The \emph{lower box-counting measure of $E$ induced by $\zeta$} is defined by
\[
  \boxm^\zeta(E)=\inf\left\{\sum\nolimits_n
  \boxm^\zeta_0(E_n):\text{$\{E_n\}$ is a countable cover of $E$}\right\}.
\]
\end{defn}
It is easy to check that $\boxm^\zeta_0$ is a Borel-regular metric pre-measure
(though is does not have to be subadditive)
and thus $\boxm^\zeta$ is a Borel-regular outer measure.
It is clear that if $\zeta$ is a gauge, then
$\hm^\zeta\leq\boxm^\zeta$.

\subsection*{Lower packing dimension}

The \emph{lower packing dimension} of $E$ is denoted and defined by
\[
  \lpdim E=\sup\{s:\boxm^s(E)>0\}.
\]
This definition is taken from~\cite{MR827889}.
Another common equivalent definition:
letting
\begin{equation}\label{boxdim}
  \lbdim E=\liminf\limits_{r\to0+}\frac{\log N_r(E)}{\abs{\log r}}
\end{equation}
we have
\[
  \lpdim E=\inf\left\{\sup\nolimits_n
  \lbdim(E_n):\text{$\{E_n\}$ is a countable cover of $E$}\right\}.
\]
It is clear that $\hdim\leq\lpdim$.

\subsection*{Porosity}
Porosity and its variations is a well-known concept.
We will deal mostly with the so called \emph{strong porosity}:
A set $E$ in a metric space $X$ is termed \emph{strongly porous}
if for each $x\in E$ there is a sequence $y_n\to x$ in $X$ such that
$\frac{d(y_n,x)}{d(y_n,E)}\to1$.
Another, equivalent definition: for $x\in E$ and $r>0$ let
\begin{align}
  \gamma(E,x,r)&=\sup\{s:\exists y\in B(x,r)\ B(y,s)\cap E=\emptyset\},\\
  \por(E,x)&=\limsup_{r\to0}\frac{\gamma(E,x,r)}{r}.
\end{align}
The set $E$ is \emph{strongly porous at $x$} if $\por(E,x)=1$;
$E$ is strongly porous if it is strongly porous at each $x\in E$;
$E$ is termed \emph{\si strongly porous} if it is a countable union
of strongly porous sets.
References:~\cite{MR943561,MR2201041}.

\subsection*{Cross product}
The following natural notion will turn useful.
\begin{defn}
Let $X,Y$ be sets and $E\subs X$, $F\subs Y$.
The \emph{cross product} of the sets $E,F$ in $X\times Y$ is the set
$E\bowtie F\subs X\times Y$ defined as follows.
\[
  E\bowtie F=(E\times Y)\cup (X\times F)=
  \{(x,y):x\in E \vee y\in F\}.
\]
Given $d\in\nset$, the \emph{cross power} of $E\subs X$
is the set $E^{\cross d}\subs X^d$ defined by
\[
  E^{\cross d}=
  E\bowtie\dots\bowtie E=
  \{(x_0,x_1,\dots,x_{d-1})\in X^d:\exists i<d\text{ such that } x_i\in E\}.
\]
\end{defn}
It is clear that the cross product depends on the ambient
sets $X,Y$, and likewise the cross power. In the sequel
occurrences of a cross product and cross power, $X$ and $Y$
are either $\Rset$ or $\Rset^d$.

\subsection*{Typical functions}
A set $E$ in a metric space is called \emph{meager} if it is
a union of countably many nowhere dense sets, and \emph{comeager},
if its complement is meager. It is clear that a set
in a metric space is comeager
if and only if it contains a dense $G_\del$-set.
A property is termed \emph{typical}
if there is a comeager set $E$ such that the property holds for all
$x\in E$. It is straightforward from the definition that a conjunction
of finitely or countably many typical properties is a typical property.

We will be concerned with typical properties of continuous functions
on a locally compact set $\Omega\subs\Rset^d$.
Recall that the space $C(\Omega)$ of continuous functions on $\Omega$
is provided with the compact-open topology whose
subbase consists of sets of the form
$\{f\in C(\Omega):f(K)\subs U\}$,
where $K\subs\Omega$ is compact and $U\subs\Rset$ is open.

In case $\Omega$ is compact, the compact-open topology coincides with
the topology of uniform convergence that is induced by the supremum norm
$\Abs f=\sup_{x\in\Omega}\abs{f(x)}$ and that the metric
on $C(\Omega)$ defined by $d(f,g)=\Abs{f-g}$ is in this case complete.
Thus, by the Baire category theorem, comeager sets in $C(\Omega)$ are dense,
and therefore meager sets are small.

In case $\Omega$ is locally compact but not compact,
the compact-open topology on $C(\Omega)$ is
not induced by the supremum norm, it is however still metrizable and complete.
General references: \cite{MR1039321} and \cite{MR953314}.
And just like in case of $\Omega$ compact, comeager sets in $C(\Omega)$ are dense,
and therefore meager sets are small.

It is common to phrase statements about typical properties as follows:
``For a typical function $f$ we have\dots''. This has to be understood
as follows: ``There is a comeager set $E\subs C(\Omega)$
such that for each $f\in E$ we have\dots''.

A subset of $\Rset^d$ with a nonempty interior is termed a \emph{domain}.

\section{Typical size of exceptional sets}\label{sec:size}

In this section we show that a typical continuous function has little lip
zero except on a very small set.

It will be convenient to generalize the notion of little lip function.
For a mapping $f\colon X\to Y$ between metric spaces, $x\in X$ and $r>0$ let
\begin{equation}\label{osc}
  \omega_f(x,r)=\diam f(B(x,r))
\end{equation}
denote the oscillation.
\begin{defn}[Generalized little lip function]
For a gauge $\phi$ define
\[
  \lip_\phi f(x)=\liminf_{r\to0}\dfrac{\omega_f(x,r)}{\phi(r)}.
\]
In the case when $\phi(r)=r^s$ for some fixed $s>0$ we write
$\lip_s f$ instead of $\lip_\phi f$ and call $\lip_s f$ the
\emph{little H\"older function}. Note that $\lip f=\lip_1 f$.
\end{defn}

We examine the sets of points where
$\lip_\phi f$ is zero, positive or infinite, respectively.
We use a shortcut notation for these sets:
\[
  \{\lip_\phi f=0\},\quad \{\lip_\phi f>0\},\quad \{\lip_\phi f=\infty\}.
\]
The following is our first major result.
\begin{thm}\label{main0}
Let $\Omega\subs\Rset^d$ be locally compact.
Let $\phi$ be a gauge and $\zeta$ a pseudo\-gauge.
For a typical function $f\in C(\Omega)$ there is $E\subs\Rset$
such that
\begin{enum}
\item $\lboxm^\zeta(E)=0$,
\item $E$ is \sistp,
\item if $x\in\Omega\setminus E^{\cross d}$, then $\lip_\phi f(x)=0$.
\end{enum}
\end{thm}
(Note that by Proposition~\ref{prop0} infra, if $\phi$ has low order,
then a stronger conclusion holds.)
\begin{proof}
We first prove the theorem for the case when $\Omega$ is compact.
We may suppose that $\Omega\subs\II d$.

We define a \emph{cube in $\II d$} as a cartesian product of
$d$ many equally long compact intervals in $[0,1]$.
For a cube $C$ denote by $\ell(C)$ the length of its side,
and for $\gamma>0$ write $\gamma C$ to denote the cube of side
$\ell(\gamma C)=\gamma\ell(C)$ that is concentric with $C$.

For each $n$ let $\Class_n$ be the collection of all finite families
$\class$ of disjoint closed cubes such that
there are $r_{\class}<1/n$, $\gamma_{\class}<1$ and $E_{\class}\subs\Rset$ such that
\begin{enumerate}
\item[(a)] $\ell(C)<\frac1n$ for all $C\in\class$.
\item[(b)]
  $\II d\setminus\bigcup_{C\in\class}\gamma^{\phantom{|}}_{\class}C\subs E_{\class}^{\cross d}$,
\item[(c)] $\lboxm^\zeta_{1/n}\bigl(E_{\class}\bigr)<2^{-n}$,
\item[(d)] $\forall x\in E_{\class}\
   \gamma(E_{\class},x,r_{\class})>(1-\frac1n)r_{\class}$.
\end{enumerate}
For each $\class\in\Class_n$ let
\[
  \Flass(\class)=
  \{f\in C(\Omega):\forall C\in\class\
  \diam f(C)<\tfrac1n \phi\bigl((1-\gamma_{\class})\ell(C)\bigr)\}
\]
and let
\begin{equation}\label{F_n}
  \Flass_n=\bigcup_{\class\in\Class_n}\Flass(\class).
\end{equation}
We will prove two claims and then derive the theorem.
\begin{claim}
$\Flass_n$ is open for each $n$.
\end{claim}
\begin{cproof}
It obviously suffices to show that
$\Flass(\class)$ is open for each $\class\in\Class_n$. Let $f\in\Flass(\class)$.
Since $\class$ is a finite collection,
there is $\eps>0$ such that
$\diam f(C)+\eps<\tfrac1n \phi\bigl((1-\gamma_{\class})\ell(C)\bigr)$ for all $C\in\class$.
If $g\in C(\Omega)$ and $\Abs{f-g}<\frac\eps2$, then
$\diam g(C)<\diam f(C)+\eps$ and thus
$\diam g(C)<\tfrac1n \phi\bigl((1-\gamma_{\class})\ell(C)\bigr)$, i.e.,
$g\in\Flass(\class)$.
\end{cproof}

\begin{claim}
$\Flass_n$ is dense for each $n$.
\end{claim}
\begin{cproof}
Fix $n$. Let $f\in C(\Omega)$ be arbitrary. Fix $\eps>0$. We will
find $g\in\Flass_n$ such that $\Abs{g-f}\leq\eps$.

We begin with constructing an appropriate $\class\in\Class_n$.
From the uniform continuity of $f$ there is $\del>0$ such that
$\diam A<\del$ implies $\diam f(A)<\eps$.
Pick an integer $k>n$ such that $1/k<\del$.
Let $\klass$ be the cover of $\II d$ by $k^d$ closed,
non-overlapping cubes of side length $1/k$
determined by the regular grid of hyperplanes with equations
$x_j=m/k$, $j=1,2,\dots,d$, $m=0,1,\dots,k$.

Since $\varliminf_{r\to0}\zeta(r)=0$, there is $\eta>0$ such that
\begin{align}
  \eta&<\frac{1}{2nk},           \label{eta0} \\
  \zeta(\eta)&<\frac{1}{2^n(k+1)}.  \label{eta}
\end{align}
Let $\beta=1-\frac{\eta k}{2}$ and define
\[
 \class=\{\beta K:K\in\klass,\ \beta K\cap\Omega\neq\emptyset\}.
\]
For $C\in\class$ we have $\ell(C)=\frac{\beta}{k}<\frac1k<\frac1n$, thus (a) holds.
Now choose $\gamma_{\class}$:
\[
  \gamma_{\class}=\frac{2-2k\eta}{2-k\eta}.
\]
If $C=\beta K\in\class$, then
\[
  \ell(\gamma_{\class}C)
  =\gamma_{\class}\beta\ell(K)=\gamma_{\class}\beta\tfrac1k=\tfrac1k-\eta.
\]
Thus the width of the gaps between neighboring cubes is $\eta$.
It follows that letting
\begin{equation}\label{Eclass}
  E_{\class}=\Bigl\{x\in[0,1]:\exists m=0,1,\dots,k\ \abs{x-\frac mk}<\frac\eta2\Bigr\}
\end{equation}
we have
\begin{equation*}\label{layers}
  \II d\setminus \bigcup_{C\in\class}\gamma_{\class}C\subs E_{\class}^{\cross d},
\end{equation*}
i.e., (b) holds.
To  prove (c), note that $E_{\class}$ is covered by $k+1$ intervals of length $\eta$.
Therefore \eqref{eta} yields
\begin{equation*}\label{eta2}
  \lboxm^\zeta_{1/n}\bigl(E_{\class}\bigr)\leq N_\eta(E_{\class})\zeta(\eta)
  < (k+1)\frac{1}{2^n(k+1)}=2^{-n},
\end{equation*}
i.e., (c) holds.

Our next goal is to choose $r_{\class}$ and prove (d). Let
\[
  r_{\class}=\frac12\left(\frac1k+\eta\right).
\]
If $x\in E_{\class}$, then there is $m$ such that $\abs{x-\frac mk}<\frac\eta2$.
Therefore $B(x,r_{\class})$ encompasses the midpoint $y$ of the
interval $I=[\frac mk+\frac\eta2,\frac {m+1}k-\frac\eta2]$ that is disjoint with
$E_{\class}$. Hence, letting $s=\frac12\left(\frac1k-\eta\right)$,
we have $\gamma(E_{\class},x,r_{\class})\geq s$.
It follows, with the aid of \eqref{eta0}, that
\[
  \frac{\gamma(E_{\class},x,r_{\class})}{r_{\class}}
  \geq \frac{s}{r_{\class}}
  \geq\frac{\frac1k-\eta}{\frac1k+\eta}=\frac{1-k\eta}{1+k\eta}
  >\frac{1-\frac1{2n}}{1+\frac1{2n}}=\frac{2n-1}{2n+1}
  >1-\frac1n
\]
as required for (d).
We proved that $\class\in\Class_n$.

The next goal is to construct $g\in\Flass(\class)$ such that $\Abs{g-f}\leq\eps$.
For each $C\in\class$, if $C\cap\Omega\neq\emptyset$,
pick $x_C\in C\cap\Omega$ and define a function on $\bigcup\class\cap\Omega$ by
\[
  h(x)=f(x_C)-f(x) \quad\text{if } x\in C\cap\Omega,\ C\in\class.
\]
Thus defined function is clearly continuous on $\bigcup\class\cap\Omega$ and since
$\diam C=\ell(C)\leq\frac{1}{k}<\del$, we have
$\Abs h\leq\eps$.
Now apply Tietze Extension Theorem:
$h$ has an extension $H\in C(\Omega)$ such that $\Abs H=\Abs h\leq\eps$.
Define $g=H+f$. Clearly $\Abs{g-f}=\Abs H\leq\eps$.
Moreover, if $x\in C\cap\Omega$, $C\in\class$, then $g(x)=h(x)+f(x)=f(x_C)$
and in particular $g$ is constant on each $C\in\class$. Consequently
$\diam g(C)=0<\frac1n\phi\bigl((1-\gamma_{\class})\ell(C)\bigr)$. It follows that
$g\in\Flass(\class)$.
\end{cproof}

The two claims prove that the set $\Flass=\bigcap_{n\in\Nset}\Flass_n$ is
a dense $G_\del$-set in $C(\Omega)$.

Let $f\in\Flass$.
By the assumption, there is a sequence $\class_n\in\Class_n$, $n\in\Nset$, such that
$f\in\Flass(\class_n)$.
Let
\begin{equation}\label{defE1}
  E=\bigcup_{n\in\Nset}\bigcap_{m\geq n} E_{\class_m}.
\end{equation}
The theorem will be proved if we show that $E$ has all of the
properties (i), (ii), (iii).

Since $\class_n\in\Class_n$, it follows from (c) that for all $n$ and $i\geq n$
\[
  \lboxm^\zeta_{1/i}\Bigl(\bigcap_{m\geq n} E_{\class_m}\Bigr)\leq
  \lboxm^\zeta_{1/i}\bigl(E_{\class_i}\bigr)<2^{-i}.
\]
Therefore $\lboxm^\zeta_0\Bigl(\bigcap_{m\geq n}E_{\class_m} \Bigr)=0$
which yields (i).

To prove (ii), we show that for each $n$ the set $F_n=\bigcap_{m\geq n}E_{\class_m}$
is \stp. Let $\rho_n=r_{\class_n}$. By (d), for each $x\in F_n$

\begin{align*}
  \por(F_n,x)
  &\geq\limsup_{m\to\infty}\frac{\gamma(F_n,x,\rho_m)}{\rho_m}\\
  &\geq\limsup_{m\to\infty}\frac{\gamma(E_{\class_m},x,\rho_m)}{\rho_m}
  \geq\limsup_{m\to\infty}1-\frac1m=1
\end{align*}
and (ii) is proved.

For each $m$ define $G_m=\bigcup_{C\in\class_m}\gamma_{\class_m}C$ and let
\[
  G=\bigcap_{n\in\Nset}\bigcup_{m\geq n} G_m.
\]
It is an easy exercise to derive from (b) that
$\Omega\setminus G\subs E^{\cross d}$. Thus for the proof of (iii) we only need
to show that $\lip_\phi f$ vanishes on $G$.

If $x\in G$, then there is an infinite set $J\subs\Nset$ such that
$x\in G_n$ for all $n\in J$.
Therefore there are, for all $n\in J$, cubes $C_n\in\class_n$ such that
$x\in \gamma_{\class_n}C_n$.
For $n\in J$ let
\begin{equation}\label{esen}
  s_n=(1-\gamma_{\class_n})\ell(C_n).
\end{equation}
Since $x\in\gamma_{\class_n}C_n$, the ball $B(x,s_n)$ is contained in
$C_n$. It follows that
\[
  \omega_f(x,s_n)=\diam f(B(x,s_n))\leq\diam f(C_n).
\]
Since $f\in\Flass(\class_n)$, the definition of $s_n$ yields
\begin{equation}\label{esen2}
  \omega_f(x,s_n)<\tfrac1n \phi\bigl((1-\gamma_{\class})\ell(C_n)\bigr)
  =\tfrac1n\phi(s_n)
\end{equation}
%
and consequently
\begin{equation*}
  \lip_\phi f(x)\leq\varlimsup_{n\in J}\frac{\omega_f(x,s_n)}{\phi(s_n)}
  \leq\lim_{n\in J}\tfrac1n=0.
\end{equation*}
This finishes the proof for the case of compact $\Omega$.

Now suppose that $\Omega$ is locally compact but not compact.
Let $\{U_n:n\in\Nset\}$ be a countable cover of $\Omega$ by
sets that are open in $\Omega$ such that $K_n=\overline U_n$
are compact and contained in $\Omega$.
By the above there is, for each $n$, a sequence of dense open sets
$\Flass_{n,j}$ in $C(K_n)$ such that if
$f\in \bigcap_{j\in\Nset}\Flass_{n,j}$,
then there is a set $E$ satisfying (i)--(iii). Let
\[
  \GG_{n,j}=\{f\in C(\Omega):f\restriction K_n\in\Flass_{n,j}\}.
\]
It is easy to show that $\GG_{n,j}$ are open and dense
in $C(\Omega)$.
Therefore $\GG=\bigcap_{n,j\in\Nset}\GG_{n,j}$
is a dense $G_\del$-set in $C(\Omega)$.

Let $f\in\GG$. Denote by $f_n=f\restriction K_n$.
For each $n$ let $E_n$ be the \si strongly porous set such that
$\lboxm^\zeta(E_n)=0$ and $\lip_\phi f_n=0$ outside $E_n^{\cross d}$
whose existence is guaranteed by the previous part of the proof.
Put $E=\bigcup_{n\in\Nset}E_n$. It is clear that $E$ is \si strongly porous
and $\lboxm^\zeta(E)=0$.

It remains to show that $\lip_\phi f=0$ outside $E^{\cross d}$.
For each $x\in\Omega$ there is $n$ such that $x\in U_n$ and since $U_n$ is open,
it follows that $\lip_\phi f(x)=\lip_\phi f_n(x)$.
Therefore
\[
  \{x\in\Omega:\lip_\phi f(x)>0\}
  \subs\bigcup_{n\in\Nset}\{x\in U_n:\lip_\phi f_n(x)>0\}
  \subs\bigcup_{n\in\Nset}E_n^{\cross d}
  \subs E^{\cross d}, \qedhere
\]
\end{proof}

\begin{thm}\label{main1}
Let $\Omega\subs\Rset^d$ be locally compact.
Let $\phi$ be a gauge and let $\gauge$ be a gauge such that
$\gauge\not\preccurlyeq d-1$.
For a typical function $f\in C(\Omega)$ there is a set $A\subs\Omega$
such that
\begin{enum}
\item $\lboxm^\gauge(A)=0$,
\item $A$ is \si strongly porous,
\item if $x\notin A$, then $\lip_\phi f(x)=0$.
\end{enum}
\end{thm}
\begin{proof}
This follows easily from the above theorem:
let $\zeta(r)=\gauge(r\sqrt d)/r^{d-1}$.
By assumption, $\varliminf_{r\to0}\zeta(r)=0$, so $\zeta$ is a pseudogauge.
Let $E\subs\Rset$ be the set guaranteed by Theorem~\ref{main0}, i.e.,
$E$ is \si strongly porous, $\lboxm^\zeta(E)=0$ and $\lip_\phi f\equiv0$
outside $E^{\cross d}$.
Let $A=E^{\cross d}\cap\Omega$.
It is clear that $A$ satisfies (iii). We need to prove (i) and (ii).

(i) We first show that, for each $B\subs\Rset$,
\begin{equation}\label{lboxm0}
  \lboxm_0^\gauge(B\times[0,1]^{d-1})\leq\lboxm_0^\zeta(B).
\end{equation}
Let $r_n\downto0$ be a sequence such that
$\lim_{n\to\infty} N_{r_n}(B)\zeta(r_n)=\lboxm_0^\zeta(B)$.
Since $N_r([0,1])\leq1+1/r$ for all $r$, letting $s_n=r_n\sqrt d$ we have
\[
  N_{s_n}(B\times[0,1]^{d-1})\leq N_{r_n}(B)(1+1/r_n)^{d-1}
\]
and thus
\begin{align*}
  \lboxm_0^\gauge(B\times[0,1]^{d-1})
  &\leq
  \varliminf N_{s_n}(B\times[0,1]^{d-1})\gauge(s_n)\\
  &\leq
  \varliminf N_{r_n}(B)\gauge(s_n)(1+1/r_n)^{d-1}\\
  &\leq
  \lim_{n\to\infty} N_{r_n}(B)(1+r_n)^{d-1}\zeta(r_n)
  \leq
  \lboxm_0^\zeta(B)
\end{align*}
and \eqref{lboxm0} follows.
By the definition,
\begin{equation}\label{sp0}
  E^{\cross d}=\bigcup_{i+j=d-1}\Rset^i\times E\times\Rset^j.
\end{equation}
Since $\lboxm^\gauge$ is \si additive, it is thus enough to show that
$\lboxm^\gauge(E\times[0,1]^{d-1})=0$ which in turn trivially follows
from \eqref{lboxm0}.

(ii)
Using \eqref{sp0} it is enough to notice the obvious fact that if $B\subs\Rset$
is \si strongly porous in $\Rset$, then $B\times[0,1]^{d-1}$ is \si strongly porous
in $\Rset^d$.
\end{proof}
\begin{thm}\label{main2}
Let $\phi$ be a gauge such that $\phi\succcurlyeq1$.
Let $\Omega\subs\Rset^d$ be a locally compact domain.
For a typical function $f\in C(\Omega)$,
$\hm^{d-1}(\{\lip_\phi f=\infty\})$ is not \si finite.
\end{thm}
This improves~\cite[Remark 4.4]{MR2213746} that only claims that
$\{\lip f=\infty\}$ is typically uncountable.
\begin{proof}
This follows easily from the following~\cite[Theorem 1.2]{MR2213746}
of Balogh and Cs\"ornyei.

\smallskip
\emph{Let $G\subs\Rset^d$ be a domain and let $f\colon G\to\Rset$ be continuous.
Assume that $\lip f<\infty$
on $G\setminus E$, where the exceptional set $E$ has \si finite
$(d{-}1)$-dimensional Hausdorff measure. Assume also that $\lip f\in L^p$
for some $p\in[1,\infty]$. Then $f\in W^{1,p}$ and if $p>d$, then
$\lip f=\Lip f=\abs{\nabla f}$ a.e.}
\smallskip

Let $G$ be any nonempty connected component of the interior of $\Omega$.
Let $\mathscr{U}\subs C(\Omega)$ be the
family of functions that are not constant on $G$. It is clearly open and dense.
Let $\psi(r)=r^d$ and
apply Theorem~\ref{main1}: There is a dense $G_\del$-set $\Flass$ such that
if $f\in\Flass$, then $\boxm^d(\{\lip f>0\})=0$ and,
in particular, $\lip f=0$ almost everywhere.
So if $f\in\mathscr{U}\cap\Flass$, then
$\lip f=0$ almost everywhere on $G$ and in particular
$\lip f\in L^\infty(G)$. So if $\hm^{d-1}(\{\lip f=\infty\})$
were \si finite, we would have, by the quoted theorem of Balogh and Cs\"ornyei,
$\nabla f=0$
almost everywhere on $G$ and $f$ would be constant on $G$.
\end{proof}

The above theorems give precise values of the
Hausdorff and lower packing dimensions of the exceptional sets
of a typical function.

\begin{thm}\label{coro2}
Let $\phi$ be any gauge.
Let $\Omega\subs\Rset^d$ be locally compact.
For a typical function $f\in C(\Omega)$,
\begin{enum}
\item $\lpdim \{\lip_\phi f>0\}\leq d-1$.
\end{enum}
If, moreover, $\Omega$ is a domain and $\phi\succcurlyeq 1$, then
\begin{enum}\setcounter{enumi}{1}
\item $\lpdim \{\lip_\phi f>0\}=\hdim \{\lip_\phi f>0\}=d-1$,
\item $\lpdim \{\lip_\phi f=\infty\}=\hdim \{\lip_\phi f=\infty\}=d-1$.
\end{enum}
\end{thm}
\begin{proof}
Let $E^+=\{x:\lip_\phi f(x)>0\}$ and $E^\infty=\{x:\lip_\phi f(x)=\infty\}$.
Define $\gauge(r)=e^{-\sqrt{\abs{\log r}}}r^{d-1}$.
Since $\gauge\not\preccurlyeq d-1$, Theorem~\ref{main1}
yields $\lboxm^\gauge(E^+)=0$.
Routine calculation shows that
for each $s>d-1$, however close to $d-1$, we have $\gauge(r)<r^s$
if $r$ is small enough, i.e., $\gauge\preccurlyeq s$.
Therefore $\lpdim E^+\leq d-1$.

On the other hand, if $\Omega$ has a nonempty interior and $\phi\succcurlyeq 1$, then
Theorem~\ref{main2} yields $\hdim E^\infty\geq d-1$, which is enough for (ii) and (iii),
since $E^\infty\subs E^+$ and $\hdim\leq\lpdim$.
\end{proof}
This theorem describes typical level sets $\{\lip f=0\}$, $\{\lip f>0\}$ and
$\{\lip f=\infty\}$. We, however, have no information about the level sets
$\{\lip f=\alpha\}$ for $0<\alpha<\infty$.
\begin{question}
Suppose that $\Omega\subs\Rset^d$ is a locally compact set with nonempty interior.
Is there a typical value of
$\hdim\{x:\lip f(x)=\alpha\}$ where $0<\alpha<\infty$?
What about $\hdim\{x:0<\lip f(x)<\infty\}$?
\end{question}

Let us point out one more consequence of Theorems~\ref{main1} and \ref{main2}.
By the classical theorem of Banach~\cite{Ban31}, a typical continuous function
$f$ has
$\Lip f(x)=\infty$ for all $x\in\Omega$. Therefore $\lip f=\Lip f$ only on the
exceptional set $\{\lip f=\infty\}$ and this exceptional set has by
Theorem~\ref{main1} a small size:

\begin{coro}
Let $\Omega\subs\Rset^d$ be locally compact.
For a typical function $f\in C(\Omega)$
\[
  \lpdim\{\lip f=\Lip f\}\leq d-1.
\]
\end{coro}

\section{H\"older spectrum}\label{sec:spectrum}
For $f\in C(\Omega)$ and $x\in\Omega$ we define the
\emph{little H\"older index} of $f$ at $x$ by
$$
  h_f(x)=\sup\{s:\lip_s f(x)=0\}=\inf\{s:\lip_s f(x)=\infty\},
$$
or, equivalently, $h_f(x)=\varlimsup\limits_{r\to0}\frac{\log \omega_f(x,r)}{\log r}$.

In this section we investigate properties of typical level sets
of little H\"older index.
First we note that if $\phi\not\succcurlyeq 1$, then $\lip_\phi f$ is
typically zero at \emph{every} point,
in contrast with Theorem~\ref{coro2}(ii) and (iii):
\begin{thm}\label{prop0}
Let $\Omega\subs\Rset^d$ be locally compact.
Let $\phi$ be a gauge. If $\phi\not\succcurlyeq 1$, then
for a typical function $f\in C(\Omega)$ we have $\lip_\phi f\equiv 0$ on $\Omega$.
\end{thm}
\begin{proof}
Suppose first that $\Omega$ is compact.
Since $\phi\not\succcurlyeq 1$ means $\varlimsup\frac{\phi(r)}{r}=\infty$,
there is a sequence $r_n\downto0$ such that $\phi(r_n)>n^2r_n$.
For each $n\in\Nset$ and $f\in C(\Omega)$ define
\[
  S_n(f)=\sup\{\abs{f(x)-f(y)}:\abs{x-y}\leq r_n\}
\]
and
\[
  \Flass_n=\{f\in C(\Omega):S_n(f)<\tfrac1n \phi(r_n)\}.
\]
Let $f\in\Flass_n$ and $\del<\frac12(\tfrac1n \phi(r_n)-S_n(f))$.
If $g\in C(\Omega)$ and $\Abs{f-g}<\del$, then
$S_n(g)<\del+S_n(f)+\del<\frac{\phi(r_n)}{n}$
and hence $g\in\Flass_n$.
We showed that $\Flass_n$ is open.

Set
\[
  \Flass_n^*=\bigcup_{m\geq n} \Flass_m,\qquad
  \Flass=\bigcap_{n\in\Nset}\Flass_n^*.
\]
Each $\Flass_n^*$ is clearly open. If $f$ is a Lipschitz function with Lipschitz
constant $L$, $m\geq\max(n,L)$ and $\abs{x-y}\leq r_m$, then
\begin{equation}\label{lip2}
  \abs{f(x)-f(y)}\leq L\abs{x-y}\leq mr_m<\tfrac1m \phi(r_m).
\end{equation}
%
and thus $S_m(f)<\tfrac1m \phi(r_m)$, i.e., $f\in\Flass_n^*$.
Since Lipschitz functions are dense in $C(\Omega)$,
we conclude that $\Flass_n^*$ is also dense in $C(\Omega)$.
It follows that $\Flass$ is a comeager set.

It remains to show that if $f\in\Flass$, then $\lip_\phi f\equiv0$ on $\Omega$.
So suppose $f\in\Flass$. Then there is an infinite set $I\subs\Nset$ such that
$f\in\bigcap_{n\in I}\Flass_n$. Thus, if $x\in\Omega$ and $n\in I$, then
$\diam f(B(x,r_n))<\frac{2\phi(r_n)}{n}$. It follows that
\[
  \lip_\phi f(x)\leq \limsup_{n\in I}\frac{\diam f(B(x,r_n))}{\phi(r_n)}
  \leq \limsup_{n\in I}\frac2n=0.
\]
This finishes the proof for the case of compact $\Omega$.

Now suppose that $\Omega$ is locally compact but not compact.
We proceed just as in the proof of Theorem~\ref{main0}.
Let $\{U_n:n\in\Nset\}$ be a countable cover of $\Omega$ by
sets that are open in $\Omega$ such that $K_n=\overline U_n$
are compact and contained in $\Omega$.
By the above there is, for each $n$, a sequence of dense open sets
$\Flass_{n,j}$ in $C(K_n)$ such that if
$f\in \bigcap_{j\in\Nset}\Flass_{n,j}$,
then $\lip_\phi f\equiv0$ on $K_n$. Let
\[
  \GG_{n,j}=\{f\in C(\Omega):f\restriction K_n\in\Flass_{n,j}\}.
\]
It is easy to show that $\GG_{n,j}$ are open and dense
in $C(\Omega)$.
Therefore $\GG=\bigcap_{n,j\in\Nset}\GG_{n,j}$
is a dense $G_\del$-set in $C(\Omega)$.

Let $f\in\mathscr G$. Denote by $f_n=f\restriction K_n$.
Clearly $\lip_\phi f_n\equiv0$ on $K_n$ for each $n$.
Let $x\in\Omega$. There is $n$ such that $x\in U_n$ and since $U_n$ is open,
it follows that $\lip_\phi f(x)=\lip_\phi f_n(x)=0$.
We proved that $\lip_\phi f(x)=0$ for each $x\in\Omega$, as required.
\end{proof}
This theorem, together with Theorem~\ref{coro2}, yield the following
information about typical $\lip_s f$.
\begin{prop}\label{coro0}
Let $\Omega\subs\Rset^d$ be locally compact.
For a typical function $f\in C(\Omega)$
\begin{enum}
\item $\lip_s f\equiv 0$ on $\Omega$ for each $s<1$.
\item There is a set $F\subs\Omega$ such that $\lpdim F\leq d-1$
and $\lip_s f\equiv 0$ on $\Omega\setminus F$ for each $s>0$.
\item If $\Omega$ is a domain, then $\{\lip_s f>0\}$
does not have \si finite measure $\hm^{d-1}$ for each $s\geq 1$.
\end{enum}
\end{prop}
\begin{proof}
(i) Let $\phi(r)=r^{1-r}$ and let $s<1$. Since $\phi\prec 1$, Proposition~\ref{prop0}
yields $\lip_\phi f\equiv 0$ on $\Omega$. On the other hand,
$s\prec\phi$ and thus $\lip_s f\equiv 0$ on $\Omega$.

(ii) Apply Theorem~\ref{coro2} with $\phi=r^{1/r}$.
For each $s>0$, however large, we have $\phi\prec s$.
Therefore if $\lip_\phi(x)=0$, then also $\lip_s f(x)=0$ for all $s>0$.
Thus it is enough to let $F=\{x\in\Omega:\lip_\phi f(x)>0\}$.

(iii) is a restatement of Theorem~\ref{main2}.
\end{proof}
We can now collect information on typical level sets of $h_f$.
\begin{prop}
Let $\Omega\subs\Rset^d$ be locally compact. Let
\begin{alignat*}{2}
  &E^f_\alpha&&=\{x\in\Omega:h_f(x)=\alpha\},  \\
  &E^f_{\leq\alpha}&&=\{x\in\Omega:h_f(x)\leq \alpha\},\quad\alpha\in[0,\infty].
\end{alignat*}
For a typical function $f\in C(\Omega)$
\begin{enum}
\item $E^f_\alpha=E^f_{\leq\alpha}=\emptyset$ if $\alpha<1$,
\item $E^f_1=E^f_{\leq1}$ and $\lpdim E^f_1\leq d-1$,
\item $\lpdim E^f_{\leq\alpha}\leq d-1$ if $1<\alpha<\infty$.
\item $\lpdim(\Omega\setminus E^f_\infty)\leq d-1$.
\end{enum}
If, moreover, $\Omega$ is a domain, then
\begin{enum}\setcounter{enumi}{4}
\item $\hdim E^f_1=\lpdim E^f_1=d-1$,
\item $\hdim E^f_{\leq\alpha}=\lpdim E^f_{\leq\alpha}=d-1$ if $1<\alpha<\infty$.
\item $\hdim(\Omega\setminus E^f_\infty)=\lpdim(\Omega\setminus E^f_\infty)=d-1$.
\end{enum}
\end{prop}
We, however, do not know anything about the sets $E^f_\alpha$ for $1<\alpha<\infty$.
\begin{question}
Suppose that $\Omega\subs\Rset^d$ is a locally compact domain.
Is there a typical value of $\hdim E^f_\alpha$ where $1<\alpha<\infty$?
\end{question}

\section{Level sets}\label{sec:level}

We now apply information on typical behavior of $\lip_\phi f$ we have gathered
to the study of level sets of a typical function.

Krishna M.~Garg~\cite{MR143855} pioneered in this research in 1963
and since then level sets have been extensively studied
by many. To name a few, Bruckner and Garg provided in~\cite{MR476939} a rather
complete topological description of level sets of a typical function;
Bernd Kirchheim~\cite{MR1260171} proved that all level sets of a typical function
$f\in C([0,1]^d)$ have Hausdorff dimension $d-1$; this result
was recently extended to compact metric spaces by
Balka, Buczolich and Elekes in the papers~\cite{MR3000710,MR3318168}.
Brian S.~Thomson~\cite{MR807990} focused on porosity and proved that all level sets
of a typical function $f\in C([0,1])$ are \si strongly porous.

In this section we show that a typical function admits a small set that
contains almost all level sets.

\begin{thm}\label{levels1}
Let $\Omega\subs\Rset^d$ be locally compact.
Let $\xi,\gauge$ be gauges and $\gauge\not\preccurlyeq d-1$.
For a typical function $f\in C(\Omega)$ there is a set $A\subs\Omega$ such that
\begin{enum}
\item $A$ is \sistp{} and $\lboxm^\gauge(A)=0$,
\item $f^{-1}(y)\subs A$ for $\hm^\xi$-almost all $y$, i.e., there is a set $B$ such that
$\hm^\xi(B)=0$ and $f\subs A\bowtie B$.
\end{enum}
\end{thm}
\begin{proof}
Let $\phi$ be a gauge satisfying $\xi(\phi(5r))\leq r^{d+1}$ for each $r>0$.

Apply Theorem~\ref{main1}: If $f$ is a typical function,
let $A=\{\lip_\phi f>0\}$ and $B=\{\lip_\phi f=0\}$.
By Theorem~\ref{main1} $A$ is \si strongly porous and $\boxm^\gauge(A)=0$
and thus (i) holds.

We now show that $\hm^\xi(f(B))=0$.
We will suppose that $\Omega$ is bounded, because otherwise
we may cover it with countably many bounded subsets $\Omega_n$ and prove
$\hm^\xi(f(B\cap\Omega_n))=0$ for each $n$.
Thus we may suppose $\Omega\subs[0,1]^d$.

For each $\del>0$ let
\[
  \VV_\del=\{B(x,r): x\in B,\ r<\del,\ \diam f(B(x,5r))<\phi(5r)<\del\}.
\]
We will need the simplest version of the Vitali $5r$-covering lemma,
(see, e.g., \cite[Theorem 2.1]{MR1333890}): The collection $\VV_\del$
is clearly a Vitali family and therefore there is
a disjoint countable collection $\{B(x_i,r_i)\}\subs\VV_\del$ such that
$\{B(x_i,5r_i)\}$ covers $B$.
For each $i$ let $E_i=f(B(x_i,5r_i))$. It follows that
$\{E_i\}$ is a $\del$-cover of $f(B)$.

Denote by $\lambda_d$
the Lebesgue measure in $\Rset^d$ and let $\alpha_d=\lambda_d(B(0,1))$
be the volume of the unit ball in $\Rset^d$. We have
\begin{multline*}
  \hm^\xi_\del(f(B))\leq\sum_i\xi(\diam E_i)\leq
  \sum_i\xi(\phi(5r_i))\leq\sum_i r_i^{d+1}\\
  =\frac1{\alpha_d}\sum_ir_i\alpha_d r_i^d
  =\frac1{\alpha_d}\sum_ir_i\lambda_d(B(x_i,r_i))
  \leq\frac\del{\alpha_d}\sum_i\lambda_d(B(x_i,r_i))
\end{multline*}
and since the balls $B(x_i,r_i)$ are disjoint
and their union is contained in $[-\del,1+\del]^d$,
it follows that
$\sum_i\lambda_d(B(x_i,r_i)\leq(1+2\del)^d$. We conclude that
\[
  \hm^\xi_\del(f(B))\leq\frac{\del(1+2\del)^d}{\alpha_d}
\]
and letting $\del\to0$ yields $\hm^\xi(f(B))=0$.

(ii) thus follows from the obvious fact that if $y\notin f(B)$,
then $f^{-1}(y)\subs A$.
\end{proof}

\section{Typical graph}\label{sec:graph}

By a classical result of Mauldin and Williams~\cite{MR860394}, the graph of
a typical function $f\in C([0,1])$ satisfies $\hdim f=1$,
see also~\cite{MR995975}.
The results of the previous sections can be rephrased in a rather surprising
way that yields this classical result: the graph
of a typical function is contained in a cross product of small sets
and it is \sistp{}.
In what follows $f$ denotes, besides a function, also its graph.
\begin{thm}\label{cross}
Let $\Omega\subs\Rset^d$ be locally compact.
Let $\xi$ be any gauge and let $\gauge$ be a gauge such that
$\gauge\not\preccurlyeq d-1$.
For a typical function $f\in C(\Omega)$ there are set
$A\subs\Omega$ and $B\subs\Rset$ such that
\begin{enum}
\item $A$ is \sistp{} and $\lboxm^\gauge(A)=0$,
\item $\hm^\xi(B)=0$
\item $f\subs A\bowtie B$.
\end{enum}
\end{thm}
\begin{proof}
Let $A$ be the set from Theorem~\ref{levels1} and
$B=\{y:f^{-1}(y)\nsubseteq A\}$. By Theorem~\ref{levels1}(ii)
$\hm^\xi(B)=0$ and clearly $f\subs A\bowtie B$.
\end{proof}
A particular choice of $\psi$ as in the proof of Theorem~\ref{coro2}
and $\xi(r)=e^{-\sqrt{\abs{\log r}}}$ yields dimension estimates:
\begin{coro}
Let $\Omega\subs\Rset^d$ be locally compact. \label{cross2}
For a typical function $f\in C(\Omega)$ there are sets
$A\subs\Omega$, $B\subs\Rset$ such that
\begin{enum}
\item $\lpdim A\leq d-1$,
\item $\hdim B=0$,
\item $f\subs A\bowtie B$.
\end{enum}
\end{coro}

Another interesting consequence of Theorems~\ref{main1} and~\ref{levels1}
claims that a typical graph is \sistp.
\begin{thm}\label{porous1}
Let $\Omega\subs\Rset^d$ be locally compact.
The graph of a typical function $f\in C(\Omega)$ is \sistp{}.
\end{thm}
The proof is based on an intuitive lemma that links porosity of $f$ with $\lip f$.
\begin{lem}
Let $\Omega\subs\Rset^d$ and $f\in C(\Omega)$. \label{porous0}
For each $x\in\Omega$
$$
  \por(f,(x,f(x)))\geq\frac{1}{1+\lip f(x)}.
$$
Consequently, if $\lip f(x)=0$, then $f$ is \stp{} at $(x,f(x))$.
\end{lem}
\begin{proof}
Let $r>0$. Write $\widehat x=(x,f(x))$ and $\omega(r)=\omega_f(x,r)$.
Clearly
\[
  f\restriction B(x,r)\subs B(x,r)\times[f(x)-\omega(r),f(x)+\omega(r)].
\]
Therefore the open ball centered at $(x,f(x)+\omega(r)+r)$ with radius $r$
is disjoint with $f$. It follows that $\gamma(f,\widehat x, \omega(r)+r)\geq r$
and thus
\begin{align*}
\por(f,(x,f(x)))
  &=\varlimsup_{r\to0}\frac{\gamma(f,\widehat x, \omega(r)+r)}{\omega(r)+r}
  \geq\varlimsup_{r\to0}\frac{r}{\omega(r)+r}\\
 & =\varlimsup_{r\to0}\frac{1}{\frac{\omega(r)}{r}+1}
  =\frac{1}{\lip f(x)+1}.
  \qedhere
\end{align*}
\end{proof}
\begin{proof}[Proof of Theorem~\ref{porous1}]
By Theorem~\ref{main1} there is a \sistp{} set $A\subs\Omega$
such that if $x\notin A$ then $\lip f=0$. The set $A\times\Rset$ is
clearly \sistp{} in $\Rset^{d+1}$. By the above Lemma~\ref{porous0}
the set $f\restriction(\Omega\setminus A)$ is \stp.
Therefore $f\subs f\restriction(\Omega\setminus A)\cup A\times\Rset$
is \sistp.
\end{proof}

\section{Application: Microscopic sets}\label{sec:micro}
We learned that the typical behavior of $\lip f$ implies that level sets and graphs
are typically small.
The aim of this section is to observe that
a typical continuous $f$ function has a microscopic exceptional set
$\{\lip f>0\}$ and also a microscopic graph.

The study of microscopic sets initiated in papers~\cite{Appell,MR1909968}
and it continues.
The following definition comes from~\cite{MR3259051,MR4122477}.
Denote by $\lambda_d$ the Lebesgue measure in $\Rset^d$.
A \emph{box} in $\Rset^d$ is a cartesian product of $d$ many compact intervals
in $\Rset$.

\begin{defn}\label{def:micro}
A set $X\subs\Rset^d$ is termed \emph{microscopic} if
for each $\eps>0$ there is a sequence of boxes $\langle B_n:n\in\nset\rangle$
covering $X$
such that $\lambda_d(B_n)\leq\eps^n$ for all $n\in\nset$.
\end{defn}

\subsection*{Properties of microscopic sets}
We first establish some elementary facts about microscopic sets.
Microscopic sets in $\Rset^d$ form a \si ideal and
every microscopic set is Lebesgue null. This is shown in~\cite{MR3259051} for $d=2$
and it easily to generalizes to arbitrary dimension.
However, microscopic sets are not
rotation invariant: by Theorem~\ref{mic2} below, the diagonal line in the plane
is not microscopic, while the $x$-axis is.

The following link of microscopic sets on the line with Hausdorff measures is proved in~\cite[5.8]{micro1}.

\begin{lem}[{\cite{micro1}}]\label{zetadef}
Let $\varsigma(r)=1/\abs{\log r}$
and $E\subs\Rset$. If $\hm^\varsigma(E)=0$, then $E$ is microscopic.
\end{lem}

A cross product of microscopic sets is microscopic. The trivial proof is omitted.
\begin{prop}\label{microcross}
If $E\subs\Rset^d$ and $F\subs\Rset^m$ are microscopic, then
$E\bowtie F\subs\Rset^{d+m}$ is microscopic.
\end{prop}

The following nice and simple characterization of microscopic sets in $\Rset^d$
generalizes~\cite[Theorem 10]{MR3259051}.
\begin{thm}\label{mic2}
A set $X\subs\Rset^d$ is microscopic if and only if there is a microscopic set
$E\subs\Rset$ such that $X\subs E^{\cross d}$.
\end{thm}
\begin{proof}
We only prove the forward implication, the other one follows from
Proposition~\ref{microcross}.
Let $\eps_k\to 0$ and fix $k\in\nset$. Since $X$ is microscopic,
there is a cover $\{C_{k,n}:n\in\nset\}$ by boxes such that
$\lambda_d(C_{k,n})\leq\eps_k^{dn}$ for all $n$.
Therefore for each $n$ there is $j(k,n)$ such that the orthogonal projection
$P_{j(k,n)}(C_{k,n})$ of $C_{k,n}$ onto the $j(k,n)$-th axis satisfies
$\diam P_{j(k,n)}(C_{k,n})\leq \sqrt[d]{\eps_k^{dn}}=\eps_k^n$.
Write
\[
  E=\bigcap_{k\in\nset}\bigcup_{n\in\nset}P_{j(k,n)}(C_{k,n}).
\]
It is clear that $E$ is a microscopic set in $\Rset$ and it takes a
straightforward calculation to prove that $X\subs E^{\cross d}$.
\end{proof}

To illustrate the size
of microscopic sets in $\Rset^d$ we show that they have Hausdorff dimension at most
$d-1$ and this estimate is sharp. Moreover, there is no converse.
\begin{thm}\label{mic1}
\begin{enum}
\item If $X\subs\Rset^d$ is microscopic, then $\hdim X\leq d-1$.
\item There is a microscopic set $X\subs\Rset^d$ such that $\hdim X=d-1$.
\item There is a compact set $C\subs\Rset^d$ such that $\hdim C=0$ and
$C$ is not microscopic.
\end{enum}
\end{thm}
\begin{proof}
(i) By Theorem~\ref{mic2} there is a microscopic set $E\subs\Rset$ such that
$X\subs E^{\cross d}$. It is a well-known fact that $\hdim E=0$ for
each microscopic $E\subs\Rset$, see~\cite{MR1909968}.
Therefore $\hdim E\times\Rset^{d-1}=d-1$ and consequently also
$\hdim E^{\cross d}=d-1$.

(ii) $\{0\}\times\Rset^{d-1}$ is by Theorem~\ref{mic2} microscopic.

(iii) As proved in~\cite{MR1909968}, there is a compact set $K\subs\Rset$ such that
$\hdim K=0$ but $K$ is not microscopic.
Let $C=\{(x,x,\dots,x)\in\Rset^d:x\in K\}$. Then clearly $\hdim C=\hdim K=0$.
If $C$ were microscopic, then there would be a microscopic set $E\subs\Rset$
such that $C\subs E^{\cross d}$. It is straightforward that the latter yields
$K\subs E$, so $E$ cannot be microscopic: a contradiction.
\end{proof}
\begin{thm}\label{mainmicro}
Let $\Omega\subs\Rset^n$ be locally compact and $\phi$ a gauge.
For a typical function $f\in C(\Omega)$, the exceptional set
$\{x\in\Omega:\lip_\phi f(x)>0\}$ is microscopic in $\Rset^d$.
\end{thm}
\begin{proof}
The theorem follows easily from Theorem~\ref{main0}. By the theorem
there is a set $E\subs\Omega$ such that
\begin{enum}
\item $\lboxm^\zeta(E)=0$,
\item if $x\in\Omega\setminus E^{\cross d}$, then $\lip_\phi f(x)=0$.
\end{enum}
Since $\hm^\zeta(E)\leq\lboxm^\zeta(E)$,
Lemma~\ref{zetadef} yields that $E$ is a microscopic set.
By Theorem~\ref{mic2}, $E^{\cross d}$ is also microscopic
and by (ii) $\{\lip_\phi>0\}\subs E^{\cross d}$, so $\{\lip_\phi>0\}$
is microscopic as well.
\end{proof}

\begin{thm}\label{mainmicro2}
Let $\Omega\subs\Rset^d$ be locally compact.
For a typical function $f\in C(\Omega)$
there is a set $A\subs\Omega$ such that
\begin{enum}
\item $A$ is microscopic in $\Rset^d$,
\item $\{y:f^{-1}(y)\nsubseteq A\}$ is microscopic in $\Rset$.
\end{enum}
\end{thm}
\begin{proof}
This is proved the same way as Theorem~\ref{levels1}, except that we use $\zeta$
(as defined by Lemma~\ref{zetadef}) in place of $\xi$ and Theorem~\ref{mainmicro}
in place of Theorem~\ref{main1} and at the end Lemma~\ref{zetadef}.
\end{proof}

\begin{thm}\label{micro4}
Let $\Omega\subs\Rset^d$ be a locally compact set.
The graph of a typical function $f\in C(\Omega)$ is microscopic in $\Rset^{d+1}$.
\end{thm}
\begin{proof}
Let $A$ be the set from Theorem~\ref{mainmicro2}. A routine check shows that
$f(\Omega\setminus A)\subs\{y:f^{-1}(y)\nsubseteq A\}$ and thus is
by Theorem~\ref{mainmicro2}(ii) microscopic.
Therefore the cross product $A\bowtie f(\Omega\setminus A)$
is by Proposition~\ref{microcross} microscopic.
Since obviously $f\subs A\bowtie f(\Omega\setminus A)$, we are done.
\end{proof}

\section*{Acknowledgment}
The author would like to thank Thomas Z\"urcher for carefully
reading the manuscript.

\bibliographystyle{amsplain}
\bibliography{typical}
\end{document}